# *The Heart is a Dust Board:*
## Abu'l Wafa Al-Buzjani, Dissection, Construction, and the Dialog Between Art and Mathematics in Medieval Islamic Culture

## Jennifer L. Nielsen

jlnr25@mail.umkc.edu
3174 SW ST RT T
Polo, MO 64671

*The Heart is a Dust Board:*
Abu'l Wafa Al-Buzjani, Dissection, Construction, and the
Dialog Between Art and Mathematics in Medieval Islamic Culture

> *"All this sewing-together and tearing-apart bewilders me,*
> *my mind reels at the variegations of affirmation and obliteration.*
> *The heart is a dust-board, he the geometer of the heart:*
> *what a marvel of figures, numbers, realities and names*
> *he inscribes!"*
> - Jelaluddin Rumi, Persian Sufi Islamic Poet,
>   "The Lovers' Tailor's Shop" (14, p. 81).

When Jelaluddin Rumi composed these lines, it is easy to imagine that the wandering sufi was standing over the shoulders of a Baghdad geometer, breathlessly watching him perform one of the most confounding dissection-construction puzzles of Abu'l Wafa al-Buzjani in the dust-laden surface of a medeival work board. The maddening beauty of the obliteration and reconstruction of complex geometric figures would have been enough to send the mystical Rumi into rapture about the nature of the creation of the universe.

Truly a geometer of the heart if there ever was one, the Persian mathematician and astronomer Abu'l Wafa al-Buzjani (AD 940-998) was awarded by his peers the title of "mohandes" geometer – a title for the most skillful and knowledgeable professional geometer of his day (*6, p. 10)*. Today he is widely considered one of the most outstanding Islamic mathematician-scientists of the tenth century. His impressive resume includes the introduction of the concepts of tangent, secant and cosecant for the first time, the compilation of indepth mathematical manuals for use by businessmen and artisans, and a collaboration with Al-Biruni to use an eclipse of the moon to determine longitude differences between Kath and Baghdad (12, p. 93; 1 pp. 9, 92,). He was also the creator of several intriguing geometric dissection-construction puzzles which have been baffling mathematicians, enlightening artists, and rescuing lucky minds from boredom for over one thousand years.

**Puzzle 1: The Tricky Triangle**

*"Draw three identical triangles, and one smaller triangle similar to them in shape, so that all four can be made into one large triangle."*
    –Abu'l Wafa Al Buzjani (2, p. 95; 11, pp. 83, 116; 13, p. 292)

This deceptively simple instruction may trick you into thinking that the larger triangle has to enclose all four triangles, but good luck trying to find such a solution. Al-Buzjani's own solution was expressed using the following slightly baffling illustration, at which a first glance may yield simple annoyance. If one is not thinking in terms of a dissection and reconstruction, the triangles do not appear to be "made into" the larger triangle at all.

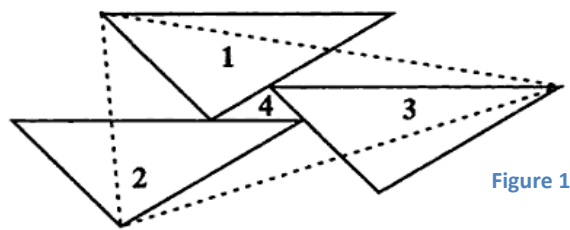

Figure 1

## Abu'l Wafa Al-Buzjani's Tricky Triangle – Steps to The Solution

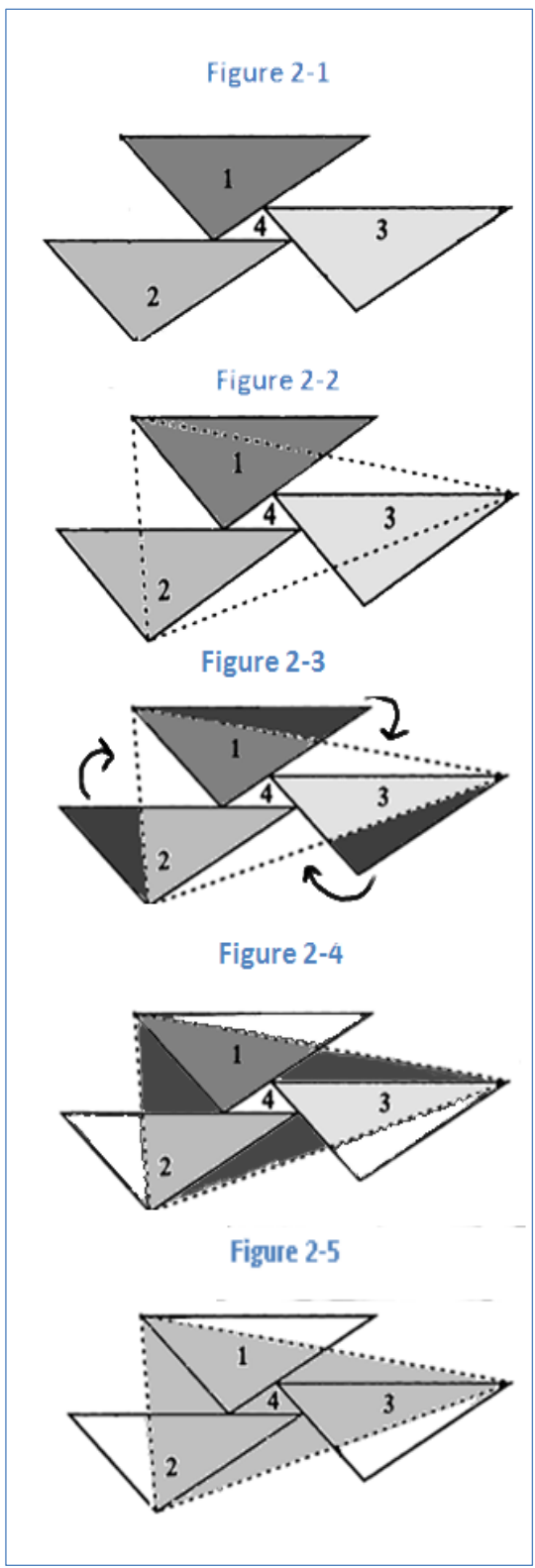

1. First place three identical triangles around a smaller similar one as shown so that the corresponding sides of the identical triangles are parallel to one another (See Figure 2-1).

2. Join three vertices of the identical triangles (one of each "type," each from a different triangle) with lines as shown to make a large triangle (Figure 2-2).

3. Now, for the "trick." Note the portions of the three identical triangles that are external to the large triangle, shown darkened (Figure 2-3). "Cut out" these three pieces and flip, rotate, or otherwise manipulate them so as to "paste" them into the "empty" portions of the large triangle.

4. See that the pieces fit exactly, as shown.

5. Now we have one large triangle and we have arrived at Al-Buzjani's solution.

**How The Tricky Triangle Works**

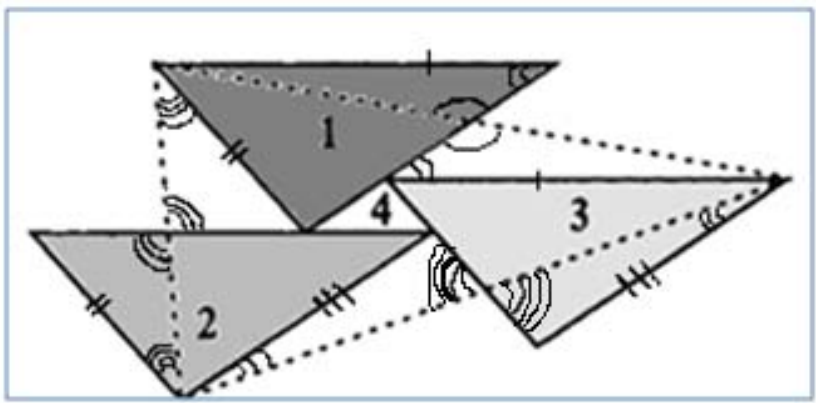

Figure 3

The "tricky triangle" can be viewed as a sliding dissection puzzle. The external portion of Triangle 2 is cut and then rotated and slid into the empty space above triangle 2 adjacent to Triangle 1. The external portion of Triangle 1 is likewise cut, rotated and slid into the empty space above Triangle 3 adjacent Triangle 1.Similarly, the external part of Triangle 3 is cut, rotated and slid into the empty space under Triangle 2 adjacent to triangle 3.

Because corresponding sides of the three identical triangles are parallel, we can see even before our operation that each "cut out" external triangle is congruent to the triangle it is pasted into by the ASA (angle-side-angle) theorem. In each case, one side of an external triangle is already known to be congruent to a side of an "empty" enclosed triangle (since these sides also happen to be corresponding sides of the original identical triangles). Furthermore, in each case, one angle of the external triangle is alternate interior to an angle in the empty triangle it is pasted into, and one angle of the external triangle is vertical to one of the angles in the empty triangle. Thus we see that what at first appears to be a conundrum reduces to a cut and paste problem utilizing simple Euclidean geometry.

**Mosaic Makers**

The tricky triangle puzzle highlights a common problem faced by tilers of mosaics in Abu'l Wafa's day: cutting apart and rearranging tile-pieces so as to assemble them into larger geometric configurations (6 p. 11). In the complex world of abstract creations in medieval Islam, a wonderful dialog was ripe for unfolding between artists and mathematicians (10, p.193). Artisans of the day often used approximations in achieving patterns, and while this was sufficient for small scale pattern working, in much larger mosaics spanning entire walls of mosques, the errors could become compounded, letting slipshod craftmanship eventually show. Thus more and more precise methods of construction and dissection such as only could be introduced by excellent mathematicians such as Al-Buzjani were necessitated. It is suspected that Al-Buzjani utilized cut and paste methods "for two purposes: to prove the correctness of certain constructions in a concrete way that could be easily understood by the artisans, and to present the constructions in such a way that the figures could be used to create new decorative patterns" some of which "became quite popular" (10, p. 193).

Buzjani wrote that "A number of geometers and artisans have made errors in the matter of…squares and their assembling. The geometers made errors because they don't have practice in applied constructing, and the artisans because they lacked knowledge of reasoning and proof." (qtd. in 6, p. 10). Al-Buzjani participated in meetings between mathematicians and artisans and was called on to give instructions in "geometric constructions of two or three dimensional ornamental patterns [and]..the application of geometry to architectural construction (qtd in 6, p. 10).

**Bad Geometry**

At one of these meetings, Buzjani explained how a common method of cutting and pasting constructions was incorrect. "Some of the artisans, [trying to create a larger square from three squares] locate one of these [three] squares in the middle and divide the next one on its diagonal and divide the third square into one isosceles right triangle and two congruent trapezoids and assemble together as it seen in the figure" (qtd in 6, p. 11).

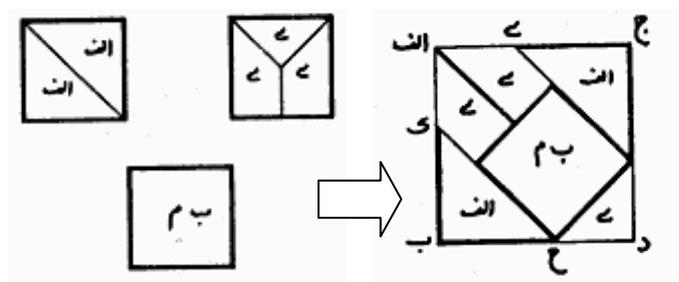
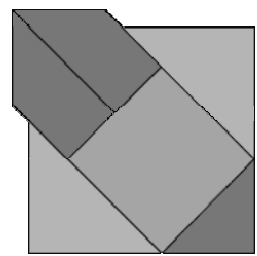

**Figure 4-1** - *Incorrect construction of a square from three unit squares. From original Persian. (Arrow added.) (6, p. 10.)*

*4-2 - Purposefully exaggerated depiction of underlying error.*

The "solution," as drawn up by the artists and in demonstration by Al-Buzjani (see figure 4-1), appears to be correct, but can be shown imprecise. It is true that the resulting shape has four right angles, and that it *appears* that each side of the larger shape is one unit plus one half of the diagonal of the unit square. But Buzjani denies that a square has truly been constructed from three squares because of this very fact, relating that "we know that each of the sides of this square is equal to the side of one of the [unit] squares plus half of its diagonal [but] it is not possible that the side of the square composed from three squares has this magnitude…[as] the diagonal of the square BG is irrational but the line HI is rational since it is equal to the side of the square BG plus half of it." It is not possible for a rational number to equal an irrational number, so the construction is "fanciful" and we "know that it is false." (7, p. 614.) An example of an easy mistake of this type is in the unit square, where a diagonal of

$\sqrt{2} \neq 1.5$, although $\sqrt{2} = 1.4$ $1 \approx 1.5$, so the mistake in construction is understandable.

**Puzzle 2: Correct Constructions of Greater Squares**

Al-Buzjani then highlighted a mathematically correct method of dividing and reassembling the squares.

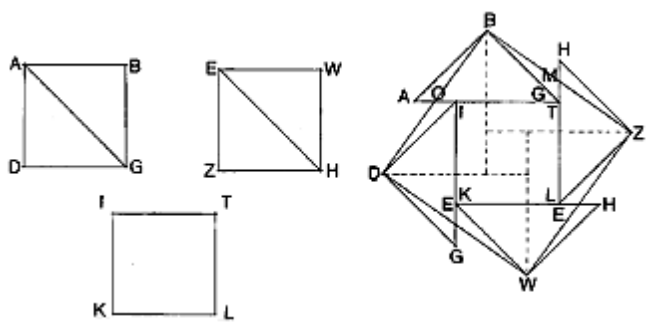

**Figure 5-1:** *The correct construction of a square from three unit squares.*

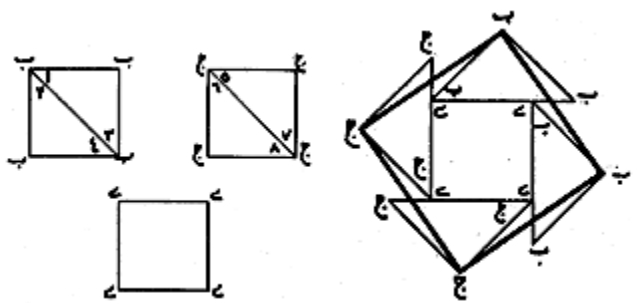

*Figure 5-2: The correct construction, as shown in the original Persian (6, p.10)*

As depicted in Figure 5-1, if we want to construct a square from three equal squares, *ABGD, EWZH, and TIKL*, we "bisect two of the squares at their diagonals by means of lines AG and EH and then transport [them] to the sides of the [third] square. Then we join the right angles of the triangles by lines *BZ, ZW, WD, DB*. On either side [of the straight line], a small triangle has now been produced from the sides of the [two big] triangles. That [empty triangle shape] is equal to the triangle which has been cut off from the big triangle. Thus triangle *BGM* is equal to triangle *MZH*, since angle *G* is half a right angle, angle *H* is half a right angle, the two opposite angles of the triangles at *M* are equal, and side *BG* is equal to side *ZH*. Therefore, the remaining sides of the triangles [*BGM, MZH*], and the triangles are equal [by the angle-side-angle theorem]." Thus, using the same basic concepts we applied in the first triangle puzzle, we may take triangle MZH and put it in the position of triangle BGM, and by the same basic argument, triangle DGE can be moved to triangle WEK, triangle ABO can be moved to triangle DIO, and triangle WHE can be moved to triangle ZLE and we have thus completed the puzzle. The figure of Al-Buzjani's dissected square became popular in Islamic ornamental arts, and can be observed today in Iran in numerous mosques including the western *iwan* of the Friday Mosque in Isfahan (10, pp. 176,177; 5).

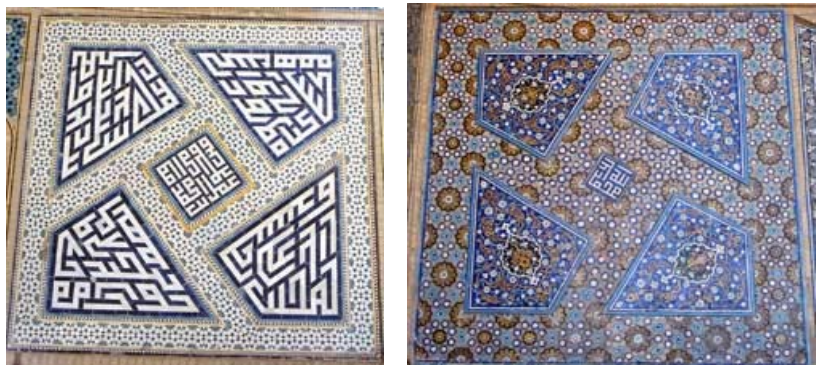

**Figure 6** - panels from the Friday Mosque in Isfahan depicting Al-Buzjani's dissected squares (10)

**A Novel Proof of the Pythagorean Theorem**

Abu'l Wafa Al-Buzjani used his dissection and construction method again to create a novel geometric proof of the Pythagorean theorem. Two unequal squares are added together to make a third square.

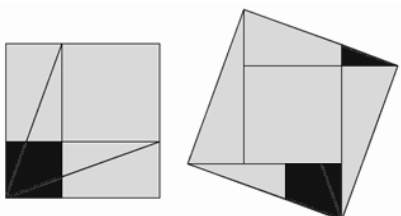 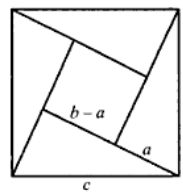

**Figure 7-1** – a small square and a large square, superimposed, are dissected and then reconstructed into a larger square, in a geometric "proof" of the Pythagorean theorem.

**Figure 7-2 –** Showing that c is the hypotenuse and that a and b are legs of the right triangle

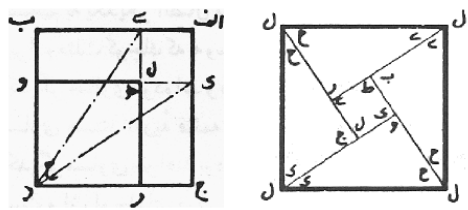

**Figure 7-3** – Drawing from original Persian text

That two equal squares are easily combined into a bigger square was known since the time of Socrates, but Abu'l Wafa's method works even if the squares are different. In Abu'l Wafa's proof, a small square is placed to share the corner and two side segments of a larger square. Then a larger triangle is built up by dissecting both the large square and the small square and adding the large square to the small square to make a bigger square. In figure 7-1, the small black square is $a^2$, the large grey square superimposed behind the small black square is $b^2$, and the larger square made from adding $a^2$ and $b^2$ is $c^2$. See Figure 7-1. In the final square, as shown in Figure 7-2, *a, b,* and *c* can be viewed as the hypotenuse and other sides of a right triangle taking up part of the square.

**Squares from More and More Squares**

Ab'l Wafa worked with many more multiples of squares in his dissections, and was interested in the general problem of assembling a "square equal in area to *n* unit squares." He distinguishes two cases in his treatise "On Assembling Squares if their Numbers is the Sum of Two Squares", being Case 1: "if the number [of unit squares] is the sum of two equal squares" (or $n = 2m^2$) and Case 2: "if the number [of unit squares] is the sum of two unequal squares" (or $n = 2ab$, where a and b are non equal, natural numbers). In case 1, his method has been described using the following process: "Cut the $2m^2$ squares in halves along a diagonal and arrange the $4m^2$ congruent triangles into a big square consisting of $m^2$ squares of area two unit squares." His method for Case 2, covering squares of different size (as we saw in the proof of the Pythagorean theorem) can be summarized as "From 2ab unit squares, compose four right-angled triangles with length a and width b. Then assemble these triangles around a square consisting of $(a-b)^2$ unit squares. If x is the side of the square obtained in this way, $x^2 = (a-b)^2 + 2ab = a^2 + b^2$." (10, p. 174).

**The Last "Puzzles"**

We have covered a small number of mentally stimulating puzzles and puzzle-proofs which were composed using the dissection-construction methods popularized by Abu'l Wafa Al-Buzjani. Many more constructions in this spirit were undertaken, generally in order to explain more complicated mathematical operations, such as solutions to cubic equations, for artisans and skilled workers who did not necessarily have prerequisite knowledge such as conics (10, p. 197). In some cases, the original, more purely mathematical arguments have been lost, while simpler and more approximate demonstrations recorded by "scribes whose training had not covered constructive geometry" (10, p. 198), and the pieces created from such arguments by such artists, are all that have survived the sands of time. This leaves us, at times, unable to answer questions about just how modern a treatment of mathematical problems the mathematicians of Al-Buzjani's time may have reached.

   We will close with a final burning question: Could Al-Buzjani and some of his contemporaries and followers have used dissection-construction methods, and other artisan-friendly treatments, to explain the concept of quasi-periodic, and even truly aperiodic, tilings, to artists?

   Aperiodic tile sets, sets of tiles allowing infinitely many distinct tilings, also known as "quasicrystals" are a phenomenon generally believed to have been first discovered by Roger Penrose in the 1970's, and are a "hot" topic of application in solid state molecular physics. The mathematics needed in a detailed investigation of this topic is beyond the scope of this paper, but a little teaser shall be left to the interested reader in the following figures.

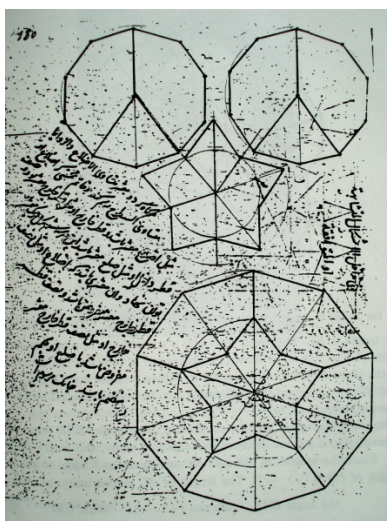
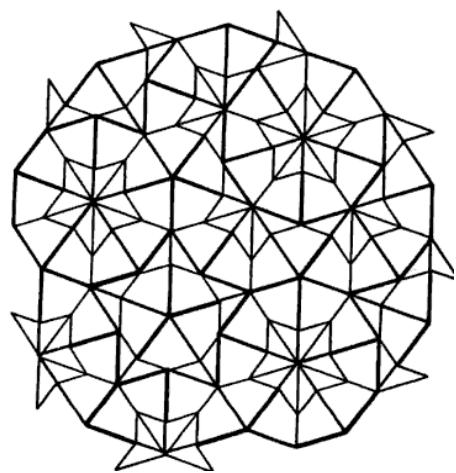

**Figure 8-1** Pentagonal Seal from Medieval Islamic Document "On Interlocking Similar Or Congruent Figures", Paris, Bibliotheque National ancien fond Persan Ms # 169 (4, p. 285)

**Figure 8-2** Aperiodic Penrose Tiling c1977, Scientific American (4, p. 286).

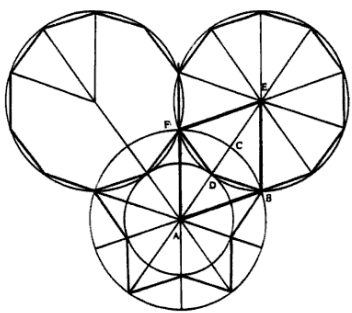

**Figure 8-3** Interlocking decagons and five pointed star (4, p. 291).

There are differences between the patterns in these pictures, to be sure, but according to one recent paper, a pentagonal pattern occurring in the tiling of the medieval Maragha mosque "is readily…obtained by a transformation of the Penrose pattern…[and] deviates from a true cartwheel Penrose tiling only in several geometric and artistic adaptations" (9, p. 85). This is enough to raise questions in the minds of several 21st century mathematicians and leading theoretical physicists (8, p. 85; 12, p.

284), and even inspire some relevant debate about the interdependence of art and science. According to mathematical historian Wasma'a K. Chorbachi,

> "[T]he true patron of the scientists who wrote these ancient manuscripts was art. It was the artisans and the architects who called for the services of science and scientists to assist them solving the design problems that they were facing. And as in the case of Islamic art in the past, science must come to the service of the arts, whether we are talking today of Islamic art, of Western art or of art generally, today more than ever before...[I]slamic tradition is so strong that, if we are in touch with the language of the present time and ground ourselves in this strong old tradition, we can arrive at an expression that is not only contemporary but could be meaningful and valid in the coming century." (8)

It can easily be argued that Al-Buzjani couldn't have known about aperiodic tiling, and it is unlikely that there were mathematical tools available at this time to fully flesh such complicated ideas out. But, as we may surmise with Rumi, who saw the creation of the universe in the swirl of numbers on the dustboard of a mathematician, Al-Buzjani was a geometer of the heart. Whether the mosaic-makers realized it or not, they were onto something special, and perhaps even in mathematics, among the intuitive ahead of their time, the heart occasionally "has its reasons of which reason knows nothing." (Blaise Pascal, 1623-1662)

May reason continually strive to catch up.

References

1. Breggren, JL. *Episodes in the Mathematics of Medieval Islam.* Springer-Verlag, 1986.
2. Danesi, Marcel. The Puzzle Instinct: The Meaning of Puzzles in Human Life. Indiana University Press, 2004.
3. Chorbachi, Wasma'a. In the Tower of Babel: Beyond symmetry in Islamic Design. Computers & Mathematics with Applications. Volume 17 , Issue 4-6 (January 1989)
4. -----  and Arthur L. Loeb. "An Islamic Pentagonal Seal." *Fivefold Symmetry.* Edited by Isvan Hargittai. World Scientific. 1992.
5. Faidherbe, Lyceé. Les Symétries du Carré. (Images 30 and 31). http://home.nordnet.fr/~ajuhel/Weyl/weyl_D4_R4.html. Accessed March 1, 2010.
6. Jaban, Slavik and Reza Sarhangi. Elementary Constructions of Persian Mosaics. R. Sarhangi and S. Jablan, *Elementary Constructions* of *Persian Mosaics*, Math Horizons, the Mathematical Association of America, September 2006, PP 10-13.
7. Imhausen, Annette and Victor J. Katz. The Mathematics of Egypt, Mesopotamia, China, India, and Islam: a Sourcebook. Princeton University Press, 2007.
8. Lu, Peter J. and Steinhardt, Paul J. "Decagonal and Quasi-Crystalline Tilings in Medieval Islamic Architecture." Science Magazine. http://www.sciencemag.org/cgi/content/full/315/5815/1106/ Accessed February 12, 2010.
9. Makovicky, Emil. "800 Year Old Pentagonal Tiling." *Fivefold Symmetry.* Edited by Isvan Hargittai. World Scientific. 1992.
10. Ozdural, Alpay. Mathematics and Arts: Connections Between Theory and Practice in the Medeival Muslim World. Historia Mathematica 27 (2000), 171–201
11. Petkovic, Miodrag. Famous Puzzles of Great Mathematicians. AMS Bookstore, 2009.
12. Suzuki, Jeff. Mathematics in Historical Context. MAA, 2009.
13. Wells, David. The Penguin Book of Curious and Interesting Puzzles. Penguin, 1992.
14. Wilson, Peter Lamborn. Sacred Drift: Essays on the Margins of Islam. City Lights Books, 1993.